\documentclass{proc-l}
\makeatletter
\AtBeginDocument{\let\@serieslogo\@empty}
\makeatother
\usepackage{amssymb}
\usepackage[hidelinks]{hyperref}
\makeatletter
\@ifundefined{subjclassname@2020}{%
  \@namedef{subjclassname@2020}{\textup{2020} Mathematics Subject Classification}%
}{}
\makeatother
\theoremstyle{plain}
\newtheorem*{thm}{Theorem}
\newtheorem{lm}{Lemma}[section]
\numberwithin{equation}{section}
\DeclareMathOperator{\Id}{Id}
\allowdisplaybreaks[4]
\begin{document}
\allowdisplaybreaks[4]
\title[Linearization for sublinear perturbations]{Linearization for unbounded perturbations of sublinear growth}
\author{Xiao Tang}
\address{School of Mathematical Sciences, Chongqing Normal University,
Chongqing 401331, China}
\email{mathtx@163.com}
\author{Jiangqiong Yu}
\address{School of Mathematical Sciences, Chongqing Normal University,
Chongqing 401331, China}
\email{jiangqiongyu168@163.com}
\thanks{Jiangqiong Yu is the corresponding author.}
\date{}
\begin{abstract}
We investigate global topological linearization with  unbounded perturbations of sublinear growth. Under a spectral gap condition and a sufficiently small Lipschitz constant, we show that the displacements of both the conjugacy and its inverse from the identity satisfy the same growth bound as the perturbation. Moreover, the conjugacy is unique within the class of maps whose displacements satisfy this bound. We obtain representations of the displacement in terms of stable and unstable orbit series that converge uniformly on any bounded set and give the inverse conjugacy by using the Banach fixed points theorem in a weighted function space.
\end{abstract}
\maketitle
\section{Introduction}
The linearization problem is a fundamental topic in the theory of
dynamical systems, as it concerns the question of whether a nonlinear
dynamical system can be transformed, by a suitable change of coordinates,
into a linear system.
Poincar\'e established local analytic linearization
of contractive differential systems near equilibria under 
nonresonance conditions; see \cite{Poincare1879}.
For one-dimensional holomorphic maps, Siegel (\cite{Siegel42}) proved analytic linearization under a
Diophantine condition on the rotation number. Brjuno (\cite{Brjuno7172})
subsequently weakened this hypothesis to the Brjuno condition.
Yoccoz (\cite{Yoccoz95}) later proved that the Brjuno condition is also
necessary for analytic linearization of quadratic polynomials near
irrationally indifferent fixed points.
For further results on analytic linearization, see
\cite{BC06,Herman87,PM92} and the references therein.
In the smooth setting, Sternberg (\cite{Sternberg57}) established smooth
linearization for local contractions under nonresonance and
differentiability assumptions. For $C^{1,1}$ local diffeomorphisms of
$\mathbb R^n$, Hartman (\cite{H60b}) showed that a contracting fixed point
admits local $C^1$ linearization without additional nonresonance
conditions. For further results on smooth linearization, see
\cite{B78,E01,RS04} and the references therein.
\par
The classical local Hartman--Grobman theorem, originating in the
work of Grobman \cite{G59,G62} and Hartman \cite{H60}, requires hyperbolicity but imposes no additional nonresonance or spectral gap conditions and requires only $C^1$ regularity.
It  states that a $C^1$ diffeomorphism admits local
topological linearization near a hyperbolic fixed point. More precisely, let $U\subset\mathbb R^n$ be an
open neighborhood of the origin and let $F:U\to\mathbb R^n$ be a $C^1$
local diffeomorphism with $F(0)=0$. If $A=DF(0)$ is hyperbolic, there is
a local homeomorphism $H$ such that
\[
H\circ F=A\circ H,
\]
wherever the compositions are defined. An analogous statement holds for
flows near hyperbolic equilibria.
\par
The Hartman--Grobman theorem has subsequently been extended in several
directions. Palis (\cite{Palis68}) treated hyperbolic fixed points in Banach
spaces, and Lu (\cite{Lu91}) obtained local linearization for scalar
reaction--diffusion equations. Results for the Cahn--Hilliard and
phase-field equations, semilinear hyperbolic evolution equations, and
nonautonomous evolution equations extend the range of infinite-dimensional
systems covered by \cite{BL94,HP16,PR16}. 
\par
Although the classical Hartman--Grobman theorem describes the dynamics
locally near a hyperbolic fixed point, it does not in general provide a
conjugacy on the entire phase space. This motivates the study of global
topological linearization. Pugh (\cite{Pugh69}) established a global
Banach-space result for bounded perturbations with sufficiently small
Lipschitz constants. Palmer (\cite{Palmer73}) subsequently developed a nonautonomous
version for differential equations admitting an exponential dichotomy. Later,
a variety of extensions in nonautonomous systems, including nonuniform hyperbolicity case, are given in \cite{AW06,BDP21,BV06,CMR07,KP91} and the references therein. 
The theory has also been extended to random dynamical systems (\cite{BV06P,ZS20}).

\par
In the bounded-perturbation setting of Pugh and Palmer, boundedness
allows one to seek a conjugacy $H=\Id+h$ with bounded continuous
displacement $h$ by a
fixed-point argument. For unbounded perturbations, one asks whether global linearization remains valid and what growth bounds can be obtained for the displacement.
Progress in this direction includes results for systems with contractive
linear parts. 
In the continuous-time setting, Lin \cite{Lin07} studied unbounded perturbations of nonautonomous differential equations under uniform contraction, while Huerta \cite{Huerta20} established linearization for such perturbations under nonuniform contraction. In the discrete-time setting, Casta\~neda and Robledo \cite{CR18} employed the method of sphere-crossing times to obtain almost topological conjugacy for uniformly asymptotically stable nonautonomous difference systems with unbounded perturbations.
When contracting and expanding directions coexist, Wu and Xia \cite{WX24} established topological linearization with partially unbounded perturbations using a two-cylinder construction. For autonomous differential equations on Banach spaces, Dragi\v{c}evi\'c, Li, Lu and Xia \cite{DLLX25} obtained global topological linearization by combining invariant foliations with a crossing-time construction. Wu, Xia and Zhou \cite{WXZ25} generalized the idea to nonautonomous difference equations on Banach spaces with unbounded Lipschitz perturbations. 
\par
In this paper, we continue to study topologically global linearization for unbounded perturbations. But sublinear growth condition is required (see definition below). We first introduce the notations.
\begingroup
Let $(X,\|\cdot\|)$ be a Banach space and let $A\in GL(X)$ be hyperbolic.
We present the estimates below for the case in which both the stable and unstable spectral subspaces of $A$ are nonzero.
Consider $F:=A+f$, where $f:X\to X$ is globally Lipschitz, may be unbounded and there are constants $0<\alpha<1$ and $C_f>0$ such that
\begin{equation}\label{eq:sublinear-growth}
     \|f(x)\| \leq C_f (1+\|x\|^\alpha) \qquad \text{for all } x\in X. 
\end{equation}
The condition \eqref{eq:sublinear-growth} is referred to as sublinear growth because it allows $\|f(x)\|$ to tend to infinity as $\|x\|\to\infty$, while ensuring that its growth rate is slower than linear when $0<\alpha<1$.
Write $X=X_s\oplus X_u$ for its stable and unstable spectral subspaces,
$P_s:X\to X_s,P_u:X\to X_u$ for the corresponding projections, and
$A_s:=A|_{X_s}$, $A_u:=A|_{X_u}$. We write $r(T)$ for the spectral radius
of a bounded linear operator $T$. 
\begin{thm}
Let $A\in GL(X)$ be hyperbolic with $X_s\ne\{0\}$ and $X_u\ne\{0\}$. Assume that $f$ is globally Lipschitz:
\begin{equation*}
     \|f(x)-f(y)\|\leq L_f\|x-y\|,
 \qquad x,y\in X,
\end{equation*}
with the Lipschitz constant $L_f>0$ sufficiently small and admits the sublinear growth property \eqref{eq:sublinear-growth} with the exponent $\alpha$ satisfying
\begin{equation}\label{eq:spectral-growth}
 \alpha<\min\left\{
 -\frac{\log r(A_s)}{\log r(A_s^{-1})},
 -\frac{\log r(A_u^{-1})}{\log r(A_u)}\right\}.
\end{equation}
Then there is a homeomorphism $H=\Id+h:X\to X$ such that
\begin{equation}\label{eq:semiconjugacy}
 H\circ F=A\circ H.
\end{equation}
Its displacement $h$ has the representations below:
\begin{align}
 P_sh(x)&=-\sum_{j=1}^{\infty}A_s^{j-1}P_sf(F^{-j}x),\label{eq:h-series-s}\\
 P_uh(x)&=\sum_{j=0}^{\infty}A_u^{-j-1}P_uf(F^jx).\label{eq:h-series-u}
\end{align}
Moreover, for a constant $C_H>0$,
\begin{equation}\label{eq:two-sided-growth}
 \|H(x)-x\|+\|H^{-1}(x)-x\|
 \leq C_H(1+\|x\|^\alpha),\qquad x\in X.
\end{equation}
The map $H$ is unique among continuous maps satisfying
\eqref{eq:semiconjugacy} whose displacements have sublinear
growth with the same exponent $\alpha$.
\end{thm}
The spectral condition \eqref{eq:spectral-growth} can equivalently be written as
\[
 r(A_s)r(A_s^{-1})^\alpha<1,\qquad
 r(A_u^{-1})r(A_u)^\alpha<1,
\]
with the inequality corresponding to a zero spectral subspace omitted. These spectral radii are unchanged by
passing to an equivalent norm, so they retain their meaning throughout
the proof. We keep $\|\cdot\|$ and $L_f$ for the original norm and
its Lipschitz bound, and use $\|\cdot\|_*$ and $L_f^*$ for the adapted
norm and a Lipschitz bound in that norm. The sufficient threshold for
the original $L_f$ is given in \eqref{eq:original-threshold} below.
\endgroup
\par
\begingroup
The theorem gives quantitative control of global linearization in a
prescribed growth class. Under the spectral condition
\eqref{eq:spectral-growth} and the small Lipschitz condition, the growth
bound \eqref{eq:sublinear-growth} on the perturbation is preserved by both
the displacement of the conjugacy and that of its inverse, as stated in
\eqref{eq:two-sided-growth}. This means that both
displacements have growth at most $1+\|x\|^\alpha$, with the same exponent
as the perturbation. In particular, the conjugacy and its inverse are asymptotic to the identity in the following sense:
\[
 \frac{\|H(x)-x\|}{\|x\|}\longrightarrow0,\qquad
 \frac{\|H^{-1}(x)-x\|}{\|x\|}\longrightarrow0
 \quad\text{as }\|x\|\to\infty.
\]
The theorem also gives uniqueness
within the prescribed displacement class.
\par
In \cite[Corollary~1]{WXZ25}, Wu, Xia and Zhou established global topological linearization for globally Lipschitz perturbations of autonomous hyperbolic linear maps on Banach spaces, without a sublinear growth bound. Their results do not provide displacement estimates of the form \eqref{eq:two-sided-growth}. Under the additional growth condition \eqref{eq:sublinear-growth} and spectral condition \eqref{eq:spectral-growth}, our theorem obtains these estimates for the conjugacy and its inverse, together with uniqueness in the prescribed growth class.
Moreover, the displacement of the conjugacy from the identity is represented by the stable and unstable orbit series \eqref{eq:h-series-s}--\eqref{eq:h-series-u}, which converge uniformly on bounded sets.

Casta\~neda and Robledo \cite[Theorems~2.6--2.7]{CR18} establish almost topological equivalence for finite-dimensional, uniformly contractive nonautonomous systems under a smallness condition on the Lipschitz constant. Their notion of equivalence allows the conjugacies and their inverses to be discontinuous on sets of Lebesgue measure zero \cite[Definition~1.3]{CR18}. By contrast, under our growth and spectral assumptions, we obtain a global homeomorphism on a Banach space whose displacement and that of its inverse satisfy the common bound \eqref{eq:two-sided-growth}, while allowing stable and unstable directions to coexist. The constructions also differ: their approach yields explicit crossing-time formulas \cite[p.~2294]{CR18}, whereas ours provides convergent orbit-series representations.
\par
The paper is organized as follows. 
Section~\ref{sec:estimates}
constructs an adapted norm, bounds the growth of forward and backward
iterates, and shows that sublinear displacement bounds are preserved
under composition and under passage to an equivalent norm.
Section~\ref{sec:main-proof} constructs the conjugacy and its
inverse and proves the theorem.
\par
\endgroup
\begingroup
\section{Useful estimates for the proof}\label{sec:estimates}
We first collect four elementary estimates. The construction of the
conjugacy and the verification that it is a homeomorphism are carried out
in Section~\ref{sec:main-proof}.
\begin{lm}[A simultaneous adapted norm]\label{lem:adapted-norm}
Let $T$ be an invertible bounded operator on a  Banach space
$(Y,\|\cdot\|)$. If $p>r(T)$ and $q>r(T^{-1})$, then
\begin{equation*}
 |y|_{p,q}:=
 \max\left\{\sup_{n\geq0}p^{-n}\|T^ny\|,
                  \sup_{n\geq1}q^{-n}\|T^{-n}y\|\right\}
\end{equation*}
is an equivalent norm, and
\begin{equation}\label{eq:simultaneous-bounds}
 |Ty|_{p,q}\leq p|y|_{p,q},\qquad
 |T^{-1}y|_{p,q}\leq q|y|_{p,q}.
\end{equation}
\end{lm}
\begin{proof}
Choose $p_0,q_0$ with $r(T)<p_0<p$ and $r(T^{-1})<q_0<q$.
The spectral radius formula gives constants $K_+,K_-$ such that
\[
 \|T^n\|\leq K_+p_0^n,\qquad
 \|T^{-n}\|\leq K_-q_0^n,\qquad n\geq0.
\]
For example, the first estimate holds with constant $1$ for all sufficiently
large $n$, after which the finitely many remaining powers can be absorbed
into $K_+$. The second estimate follows in the same way. Therefore
\[
 \|y\|\leq |y|_{p,q}\leq C_T\|y\|
\]
for some $C_T\geq1$. Homogeneity and the triangle inequality follow from
the corresponding properties of each term in the suprema, and the
$n=0$ term shows definiteness.
\par
For every $n\geq1$, $\Id=T^nT^{-n}$ gives
$1\leq\|T^n\|\,\|T^{-n}\|$. Taking $n$th roots and letting $n\to\infty$ yields
$1\leq r(T)r(T^{-1})<pq$.
Thus $q^{-1}<p$ and $p^{-1}<q$. For the positive powers,
\[
 p^{-n}\|T^{n+1}y\|\leq p|y|_{p,q},\qquad n\geq0.
\]
For the negative powers, the term $n=1$ is bounded by
$q^{-1}\|y\|\leq p|y|_{p,q}$, and for $n\geq2$,
\[
 q^{-n}\|T^{-(n-1)}y\|\leq q^{-1}|y|_{p,q}
 \leq p|y|_{p,q}.
\]
This proves the first inequality in \eqref{eq:simultaneous-bounds}.
For $T^{-1}$, the $n=0$ term in the first supremum satisfies
$\|T^{-1}y\|\leq q|y|_{p,q}$. For $n\geq1$,
\[
 p^{-n}\|T^{n-1}y\|\leq p^{-1}|y|_{p,q}\leq q|y|_{p,q}.
\]
The terms in the second supremum satisfy
\[
 q^{-n}\|T^{-(n+1)}y\|\leq q|y|_{p,q},\qquad n\geq1.
\]
This proves the other inequality.
\end{proof}
\begin{lm}[Growth of forward and backward iterates]\label{lem:orbit-growth}
Let $(X,\|\cdot\|)$ be a Banach space and $A\in GL(X)$ be hyperbolic with $X_s\ne\{0\}$ and $X_u\ne\{0\}$.
Let $0<\alpha<1$ and assume that $F=A+f$ is a
homeomorphism with
\[
 \|f(x)\|\leq C(1+\|x\|^\alpha),\qquad x\in X.
\]
For any $B>\|A^{-1}\|$ and
$D>\|A\|$, there is $R\geq1$ such that, for every $j\geq0$,
\begin{align}
 R+\|F^{-j}x\|
 &\leq B^{j}(R+\|x\|),\notag\\
 R+\|F^jx\|
 &\leq D^{ j}(R+\|x\|),\label{eq:nonlinear-weights}\\
 R+\|A^{-j}x\|
 &\leq \|A^{-1}\|^{ j}(R+\|x\|),\notag\\
 R+\|A^jx\|
 &\leq \|A\|^{j}(R+\|x\|).\label{eq:linear-weights}
\end{align}
Moreover, for all $0<\alpha<1$,
\begin{equation}\label{eq:explicit-growth-comparison}
 \tfrac12(1+\|x\|^\alpha)\leq(R+\|x\|)^\alpha
 \leq R^\alpha(1+\|x\|^\alpha).
\end{equation}
\end{lm}
\begin{proof}
Since $f$ satisfies the sublinear growth property, for every $\eta>0$ there is $C_\eta\geq0$ such that
\begin{equation}\label{eq:sublinear-absorption}
 \|f(x)\|\leq\eta\|x\|+C_\eta.
\end{equation}
Indeed, choose $\tilde \eta\geq1$ such that
$C\tilde \eta^{\alpha-1}\leq\eta$. 
Since $\alpha-1<0$, for $\|x\|\geq \tilde \eta$,
we have
\(
\|x\|^{\alpha-1}\leq \tilde \eta^{\alpha-1}.
\)
Multiplying both sides by $C\|x\|$,
we obtain
\[
\begin{aligned}
C\|x\|^\alpha
=C\|x\|^{\alpha-1}\|x\|
\leq C\tilde \eta^{\alpha-1}\|x\|
\leq\eta\|x\|.
\end{aligned}
\]
Using the sublinear growth condition gives
\[
\begin{aligned}
\|f(x)\|
\leq C(1+\|x\|^\alpha)
=C+C\|x\|^\alpha
\leq C+\eta\|x\|.
\end{aligned}
\]
For
$\|x\|\leq \tilde \eta$,
by the sublinear growth condition,  we have
\begin{equation*}
  \|f(x)\|
\leq C(1+\|x\|^\alpha)
\leq C(1+\tilde \eta^\alpha) \le \eta\|x\|+C(1+\tilde \eta^\alpha).  
\end{equation*}
Thus \eqref{eq:sublinear-absorption} holds with
$C_\eta:=C(1+\tilde \eta^\alpha)$.
\par
Taking $\eta=D-\|A\|$ in \eqref{eq:sublinear-absorption}, we obtain
\begin{equation*}
 \|Fx\|\leq D\|x\|+C_+.
\end{equation*}
Indeed, $\eta>0$ because $D>\|A\|$, and the preceding estimate gives
\[
 \|Fx\|\leq\|A\|\,\|x\|+\|f(x)\|
 \leq(\|A\|+\eta)\|x\|+C_\eta
 =D\|x\|+C_\eta.
\]
Thus we may take $C_+=C_\eta$ for this choice of $\eta$.
\par
For the backward estimate, put $y:=F^{-1}x$ and hence $Ay=x-f(y)$.
With $\eta=\|A^{-1}\|^{-1}-B^{-1}>0$, we have
\[
 \begin{aligned}
 \|y\|&\leq\|A^{-1}\|\,\|x\|
                  +\|A^{-1}\|\,\|f(y)\|\\
 &\leq\|A^{-1}\|\,\|x\|
       +\|A^{-1}\|\eta\|y\|
       +\|A^{-1}\|C_\eta.
 \end{aligned}
\]
Since $1-\|A^{-1}\|\eta=\|A^{-1}\|/B$, this gives
\begin{equation}\label{eq:backward-linear-growth}
 \|F^{-1}x\|\leq B\|x\|+C_-,\qquad C_-=BC_\eta.
\end{equation}
More explicitly, moving the term involving $\|y\|$ to the left gives
\[
 (1-\|A^{-1}\|\eta)\|y\|
 \leq\|A^{-1}\|\,\|x\|+\|A^{-1}\|C_\eta.
\]
Since $1-\|A^{-1}\|\eta=\|A^{-1}\|/B>0$, division by this number yields
\[
 \begin{aligned}
 \|y\|
 &\leq\frac{\|A^{-1}\|}{1-\|A^{-1}\|\eta}\|x\|
       +\frac{\|A^{-1}\|C_\eta}{1-\|A^{-1}\|\eta}\\
 &=B\|x\|+BC_\eta.
 \end{aligned}
\]
This is \eqref{eq:backward-linear-growth}, since $y=F^{-1}x$.
Only the assumed existence of the continuous inverse is used here;
no Lipschitz bound for $F^{-1}$ is required in this estimate.
\par
Choose
\[
 R\geq\max\{1,C_+/(D-1),C_-/(B-1)\}.
\]
Then
\[
 R+\|Fx\|\leq D(R+\|x\|),\qquad
 R+\|F^{-1}x\|\leq B(R+\|x\|).
\]
For example, the first inequality follows from
\[
 \begin{aligned}
 R+\|Fx\|&\leq R+D\|x\|+C_+\\
 &\leq R+D\|x\|+(D-1)R=D(R+\|x\|),
 \end{aligned}
\]
and the second follows in exactly the same way from
$C_-\leq(B-1)R$. The case $j=0$ is equality. If the forward bound holds
for $j$, then
\[
 R+\|F^{j+1}x\|
 \leq D(R+\|F^jx\|)
 \leq D^{j+1}(R+\|x\|).
\]
The backward bound follows by the same induction, using
\[
 R+\|F^{-(j+1)}x\|
 \leq B(R+\|F^{-j}x\|)
 \leq B^{j+1}(R+\|x\|).
\]
Thus, the induction gives \eqref{eq:nonlinear-weights}.
Clearly, $R+\|A^{\pm1}x\|$ is bounded by
$\|A\|(R+\|x\|)$ and
$\|A^{-1}\|(R+\|x\|)$, respectively.
Therefore the induction gives
\[
 R+\|A^{\pm j}x\|
 \leq\|A^{\pm1}\|^{j}(R+\|x\|),\qquad j\geq0,
\]
proving
\eqref{eq:linear-weights}
\par
The left inequality in
\eqref{eq:explicit-growth-comparison} follows because
$(R+\|x\|)^\alpha$ is at least both $1$ and $\|x\|^\alpha$.
The right one follows from $(r+t)^\alpha\leq r^\alpha+t^\alpha$
and $R\geq1$.
Explicitly,
\[
 (R+\|x\|)^\alpha
 \geq\max\{R^\alpha,\|x\|^\alpha\}
 \geq\max\{1,\|x\|^\alpha\}
 \geq\tfrac12(1+\|x\|^\alpha),
\]
whereas
\[
 (R+\|x\|)^\alpha
 \leq R^\alpha+\|x\|^\alpha
 \leq R^\alpha(1+\|x\|^\alpha).
\]
This completes the proof.
\end{proof}
\begin{lm}[Growth under composition]\label{lem:composition-growth}
Let $(X,\|\cdot\|)$ be a normed space and $0<\alpha<1$.
Suppose that $u,v:X\to X$ satisfy
\[
 \|u(x)\|\leq C_u(1+\|x\|^\alpha),\qquad
 \|v(x)\|\leq C_v(1+\|x\|^\alpha).
\]
Then $(\Id+v)\circ(\Id+u)-\Id$ satisfies the same growth bound
up to a constant.
\end{lm}
\begin{proof}
The inequalities $1+\|x\|^\alpha\geq1$
and $(r+t)^\alpha\leq r^\alpha+t^\alpha$ imply
\begin{align*}
 1+\|x+u(x)\|^\alpha
 &\leq 1+\|x\|^\alpha+\|u(x)\|^\alpha\\
 &\leq 1+\|x\|^\alpha+C_u^\alpha(1+\|x\|^\alpha)^\alpha\\
 &\leq (1+C_u^\alpha)(1+\|x\|^\alpha).
\end{align*}
Consequently,
\begin{align*}
   \|\{(\Id+v)\circ(\Id+u)-\Id \}(x)\| &=  \|u(x)+v(x+u(x))\|
   \\
 &\leq\bigl(C_u+C_v(1+C_u^\alpha)\bigr)(1+\|x\|^\alpha).
\end{align*}
This completes the proof.
\end{proof}

\begin{lm}[Equivalent sublinear growth bounds]
\label{lem:equivalent-sublinear-growth}
Let $(X,\|\cdot\|)$ be a normed space,  $u:X\to X$,
and $0<\alpha<1$. Fix $R\ge 1$. Then the following
two conditions are equivalent:
\begin{enumerate}
    \item[(i)] There exists a constant $C>0$ such that
    \[
        \|u(x)\|\le C\bigl(1+\|x\|^\alpha\bigr)
        \qquad\text{for all }x\in X.
    \]

    \item[(ii)] There exists a constant $C_R>0$ such that
    \[
        \|u(x)\|\le C_R\bigl(R+\|x\|^\alpha\bigr)
        \qquad\text{for all }x\in X.
    \]
\end{enumerate}
The constants $C$ and $C_R$ both are independent of $x$.
\end{lm}

\begin{proof}
We first prove the equivalence of \textnormal{(i)} and
\textnormal{(ii)}. Since $R>1$, for every $t\ge0$ we have
\[
    1+t^\alpha
    \le R+t^\alpha
    \le R(1+t^\alpha).
\]
Taking $t=\|x\|$, the first inequality shows that
\textnormal{(i)} implies \textnormal{(ii)} with $C_R=C$,
while the second shows that \textnormal{(ii)} implies
\textnormal{(i)} with $C=RC_R$.
The proof is completed.
\end{proof}

\section{Proof of the Theorem}\label{sec:main-proof}
\begin{proof} 
By \eqref{eq:spectral-growth}, choose $\epsilon>0$ small enough so that
\begin{equation*}
\mathcal{C}_ s:=r(A_s)+\epsilon,\qquad \mathcal{C}_u:=r(A_u^{-1})+\epsilon
\end{equation*}
satisfy $0<\mathcal{C}_s<1$, $0<\mathcal{C}_u<1$, and
\begin{equation}\label{contrac.with-alpha}
 \mathcal{C}_s\bigl(r(A_s^{-1})+\epsilon\bigr)^\alpha<1,\qquad
 \mathcal{C}_u\bigl(r(A_u)+\epsilon\bigr)^\alpha<1.    
\end{equation}
Such a choice is possible by continuity in $\epsilon$ and the two strict
spectral inequalities. 
Now we construct an equivalent norm $\|\cdot\|_*$ on $X$ so that 
\begin{equation}\label{coeff.-contraction}
\|A^{-1}\|_*L_f^*<1\qquad \text{and} \qquad L_f^*\kappa<1,    
\end{equation}
where $L_f^*$ is the Lipschitz constant of $f$ under the new norm $\|\cdot\|_*$ and  
\begin{equation*}
 \kappa:=\max\left\{\frac{\|A^{-1}\|_*^\alpha}{1-\mathcal{C}_s\|A^{-1}\|_*^\alpha},
                         \frac{\mathcal{C}_u}{1-\mathcal{C}_u\|A\|_*^\alpha}\right\}.
\end{equation*}

Apply Lemma~\ref{lem:adapted-norm} to $(A_s,\mathcal{C}_s,r(A_s^{-1})+\epsilon)$ and $(A_u,r(A_u)+\epsilon,\mathcal{C}_u)$, we can define new norms $\|\cdot \|_s$ and $\|\cdot\|_u$ on stable subspace $X_s$ and unstable subspace $X_u$ respectively, both of which are equivalent to the original one, such that
\begin{align}
 \|A_s\|_s&\leq \mathcal{C}_s,\qquad
 \|A_s^{-1}\|_s\leq r(A_s^{-1})+\epsilon,\label{est.norm.linar.contrac.}\\   
 \|A_u\|_u&\leq r(A_u)+\epsilon,\qquad
 \|A_u^{-1}\|_u\leq \mathcal{C}_u.\label{est.norm.linear.expand.}
\end{align}
Set
\begin{equation*}
 \|x\|_*:=\max\{\|P_sx\|_s,\|P_ux\|_u\}.
\end{equation*}
This is an equivalent norm on $X$, so there are $m,M>0$ such that
\begin{equation}\label{eq:norm-equivalence}
 m\|x\|\leq\|x\|_*\leq M\|x\|,\qquad x\in X.
\end{equation}
Indeed, Lemma~\ref{lem:adapted-norm} gives constants $C_s,C_u\geq1$ such that
$\|z\|\leq\|z\|_s\leq C_s\|z\|$ and
$\|z\|\leq\|z\|_u\leq C_u\|z\|$, respectively. Hence
\[
 \|x\|\leq\|P_sx\|+\|P_ux\|\leq2\|x\|_*,\qquad
 \|x\|_*\leq
 \max\{C_s\|P_s\|,C_u\|P_u\|\}\|x\|.
\]
Thus, one may take $m=1/2$
and $M=\max\{C_s\|P_s\|,C_u\|P_u\|\}$, where the projection norms
in this expression are computed in the original norm.
Moreover, $\|P_sx\|_*\leq\|x\|_*$ and
$\|P_ux\|_*\leq\|x\|_*$, so projections introduce no additional
constant in the component estimates below.
Since $A$ is block diagonal for the splitting $X_s\oplus X_u$,
the maximum norm also gives
\begin{align*}
\|A\|_*&=\max\{\|A_s\|_s,\|A_u\|_u\}=\|A_u\|_u,
\\
 \|A^{-1}\|_*&=\max\{\|A_s^{-1}\|_s,\|A_u^{-1}\|_u\}=\|A_s^{-1}\|_s.
\end{align*}
It follows from \eqref{contrac.with-alpha}, \eqref{est.norm.linar.contrac.} and \eqref{est.norm.linear.expand.} that
\begin{equation}\label{eq:norm-gaps}
 \mathcal{C}_s\|A^{-1}\|_*^\alpha<1,\qquad
 \mathcal{C}_u\|A\|_*^\alpha<1.
\end{equation}
\par
The assumptions on $f$ give
\begin{equation*}
 \|f(x)-f(y)\|_*\leq L_f^*\|x-y\|_*,\qquad
 L_f^*:=\frac{M}{m}L_f.
\end{equation*}
Here the Lipschitz estimate follows from
\[
 \|f(x)-f(y)\|_*
 \leq ML_f\|x-y\|
 \leq\frac{M}{m}L_f\|x-y\|_*.
\]
Similarly,
\[
 \|f(x)\|_*
 \leq MC_f(1+\|x\|^\alpha)
 \leq MC_f\bigl(1+m^{-\alpha}\|x\|_*^\alpha\bigr).
\]
Thus $\|f(x)\|_*\leq C(1+\|x\|_*^\alpha)$, where one may take
$C=MC_f\max\{1,m^{-\alpha}\}$.

A sufficient threshold for the Lipschitz constant $L_f$ is
\begin{equation}\label{eq:original-threshold}
 \delta:=\frac{m}{M}\min\{\|A^{-1}\|_*^{-1},\kappa^{-1}\}.
\end{equation}
The threshold $\delta$ depends only on $A$, $\alpha$ and the original
norm. Assume that the Lipschitz constant $L_f$ of $f$ under the original norm is less than $\delta$; then \eqref{coeff.-contraction} holds.

Now we check that $F$ is a homeomorphism on $X$ if $L_f$ is chosen to satisfy $L_f<\delta$. Indeed,
For each $y\in X$, the map $x\mapsto A^{-1}(y-f(x))$ is a
contraction on $(X,\|\cdot\|_*)$ with constant at most $\|A^{-1}\|_*L_f^*<1$.
Its unique fixed point solves $F(x)=y$, so $F$ is bijective. Also,
\[
 \|F(x)-F(z)\|_*
 \geq(\|A^{-1}\|_*^{-1}-L_f^*)\|x-z\|_*,\qquad x,z\in X,
\]
which gives $\operatorname{Lip}_*(F^{-1})\leq \|A^{-1}\|_*/(1-\|A^{-1}\|_*L_f^*)$.
Indeed, $x-z=A^{-1}A(x-z)$ implies
\[
 \|A(x-z)\|_*\geq\|A^{-1}\|_*^{-1}\|x-z\|_*.
\]
Subtracting the Lipschitz bound for $f(x)-f(z)$ gives the preceding
lower estimate. Its coefficient is positive by the choice of $\delta$;
applying it to $x=F^{-1}y$ and $z=F^{-1}w$ proves continuity of the
inverse. The map $F$ itself is Lipschitz, hence continuous.
Thus $F$ is a homeomorphism.
\par
To simplify notation, write $f_s:=P_sf$, $f_u:=P_uf$, and use subscripts
for the components of other maps as well. 
\par
\medskip
\noindent\textit{Step 1. Choice of the growth bounds.}
Since the inequalities \eqref{eq:norm-gaps} and $L_f^*\kappa<1$ are strict,
we may choose $B>\|A^{-1}\|_*$ and
$D>\|A\|_*$, sufficiently close to $\|A^{-1}\|_*$ and $\|A\|_*$ respectively, such that
\begin{gather}
 \mathcal{C}_sB^\alpha<1,\qquad \mathcal{C}_uD^\alpha<1,\label{eq:buffered-gaps}\\
 L_f^*\kappa_{B,D}<1,\qquad
 \kappa_{B,D}:=
 \max\left\{\frac{B^\alpha}{1-\mathcal{C}_sB^\alpha},
                    \frac{\mathcal{C}_u}{1-\mathcal{C}_u
    D^\alpha}\right\}.
 \label{eq:buffered-smallness}
\end{gather}
Fix the corresponding $R$ from Lemma~\ref{lem:orbit-growth}.
In that lemma we now take the norm to be $\|\cdot\|_*$.
We claim that there is $C_{R}>0$ such that
\begin{equation}
 \|f(x)\|_*\leq C_{R} (R+\|x\|_*)^\alpha,\qquad x\in X.\label{new-subliear-growth}
\end{equation}
Indeed, 
by equivalence of the norms,
\[
\|f(x)\|_* \le M\|f(x)\|,
\qquad
\|x\| \le m^{-1}\|x\|_*.
\]
Since $\alpha>0$, the second inequality implies
\[
\|x\|^\alpha \le m^{-\alpha}\|x\|_*^\alpha.
\]
Consequently, for every $x\in X$,
\[
\begin{aligned}
\|f(x)\|_*
&\le M\|f(x)\|\\
&\le MC_f\bigl(1+\|x\|^\alpha\bigr)\\
&\le MC_f\bigl(1+m^{-\alpha}\|x\|_*^\alpha\bigr)\\
&\le MC_f\max\{1,m^{-\alpha}\}
       \bigl(1+\|x\|_*^\alpha\bigr)\\
&=:C_{f,*}\bigl(1+\|x\|_*^\alpha\bigr).
\end{aligned}
\]
It follows from \eqref{eq:explicit-growth-comparison} that
\[
 \|f(x)\|_*\leq C_{f,*}(1+\|x\|_*^\alpha)
 \leq 2C_{f,*}(R+\|x\|_*)^\alpha.
\]
The claim is proved with $C_R:=2C_{f,*}$.
\par
\medskip
\noindent\textit{Step 2. Construction of the conjugating map.}
If we seek $H=\Id+h$, equation \eqref{eq:semiconjugacy} becomes
\begin{equation}\label{reduced-equat.}
 h(Fx)-Ah(x)=-f(x).
\end{equation}
In the stable component, backward iteration suggests the first series
in \eqref{eq:h-series-s}; in the unstable component, forward iteration
suggests \eqref{eq:h-series-u}. We define $h$ by those two series and
first justify their convergence. By \eqref{est.norm.linar.contrac.}, \eqref{new-subliear-growth}, \eqref{eq:nonlinear-weights}, \eqref{eq:buffered-gaps} and \eqref{eq:explicit-growth-comparison},
\begin{align}
 \sum_{j=1}^{\infty}\|A_s^{j-1}f_s(F^{-j}x)\|_*
 &\leq C_{R}\sum_{j=1}^{\infty}\mathcal{C}_s^{j-1}(R+\|F^{-j}x\|_*)^\alpha\notag\\
 &\leq C_{R}\sum_{j=1}^{\infty}\mathcal{C}_s^{j-1}B^{\alpha j}(R+\|x\|_*)^\alpha\notag\\
 &=\frac{C_{R}B^\alpha}{1-\mathcal{C}_sB^\alpha}(R+\|x\|_*)^\alpha,\label{eq:stable-series-bound}\\
 \sum_{j=0}^{\infty}\|A_u^{-j-1}f_u(F^jx)\|_*
 &\leq C_{R}\sum_{j=0}^{\infty}\mathcal{C}_u^{j+1}D^{\alpha j}(R+\|x\|_*)^\alpha\notag\\
 &=\frac{C_{R}\mathcal{C}_u}{1-\mathcal{C}_uD^\alpha}(R+\|x\|_*)^\alpha.\label{eq:unstable-series-bound}
\end{align}
This shows that the two series \eqref{eq:h-series-s} and \eqref{eq:h-series-u} are both absolutely convergent on any bounded subset of $X$. Every summand is
continuous, so the limit $h$ is continuous.  Moreover,
\begin{equation}\label{eq:h-estimate}
 \|h(x)\|_*\leq C_{R}\kappa_{B,D}(R+\|x\|_*)^\alpha.
\end{equation}
\par
Reindexing the absolutely convergent series gives
\begin{align*}
 h_s(Fx)
 &=-\sum_{j=1}^{\infty}A_s^{j-1}f_s(F^{1-j}x)\\
 &=-f_s(x)-\sum_{j=1}^{\infty}A_s^jf_s(F^{-j}x)
 =-f_s(x)+A_sh_s(x),\\
 h_u(Fx)
 &=\sum_{j=0}^{\infty}A_u^{-j-1}f_u(F^{j+1}x)
 =\sum_{j=1}^{\infty}A_u^{-j}f_u(F^jx)\\
 &=A_uh_u(x)-f_u(x).
\end{align*}
Thus \eqref{reduced-equat.} holds and $H\circ F=A\circ H$.
At this point $H$ is continuous, but its bijectivity remains to be proved.
\par
\medskip
\noindent\textit{Step 3. Construction of the map in the reverse direction.}
We also seek $G=\Id+g$ satisfying
\begin{equation*}
 G\circ A=F\circ G,
\end{equation*}
which is equivalent to 
\[g(Ax)-Ag(x)=f(x+g(x)). \]
Consider the space of continuous maps
\[
 \mathcal C_{\alpha,R}=
 \left\{g:X\to X:
       \|g\|_{\alpha,R}:=\sup_{x\in X}\frac{\|g(x)\|_*}{(R+\|x\|_*)^\alpha}<\infty\right\}.
\]
This is a Banach space: division by the positive continuous function
$(R+\|x\|_*)^\alpha$ identifies it isometrically with the space of
bounded continuous $X$-valued functions, complete in the supremum norm.
\par
For $g\in\mathcal C_{\alpha,R}$, define the operator $\mathcal T$ by
\begin{align*}
 (\mathcal Tg)_s(x)
 &=\sum_{j=1}^{\infty}A_s^{j-1}
       f_s(A^{-j}x+g(A^{-j}x)),\\
 (\mathcal Tg)_u(x)
 &=-\sum_{j=0}^{\infty}A_u^{-j-1}
       f_u(A^jx+g(A^jx)).
\end{align*}
We claim that $\mathcal{T}$ has a unique fixed point in the space $\mathcal{C}_{\alpha,R}$. Indeed,
for any $z\in X$, using the Lipschitz property of $f$ yields
\begin{align*}
 \|f(z+g(z))\|_* &= \|f(z)+ f(z+g(z)) -f(z) \|_*
 \\
 &\leq\|f(z)\|_*+L_f^*\|g(z)\|_*
 \\
 &\leq(C_{R}+L_f^*\|g\|_{\alpha,R})(R+\|z\|_*)^\alpha.    
\end{align*}
Applying \eqref{eq:linear-weights} in each series gives
\begin{align*}
 \|(\mathcal Tg)_s(x)\|_*
 &\leq(C_{R}+L_f^*\|g\|_{\alpha,R})
       \frac{\|A^{-1}\|_*^\alpha}{1-\mathcal{C}_s\|A^{-1}\|_*^\alpha}(R+\|x\|_*)^\alpha,\\
 \|(\mathcal Tg)_u(x)\|_*
 &\leq(C_{R}+L_f^*\|g\|_{\alpha,R})
       \frac{\mathcal{C}_u}{1-\mathcal{C}_u\|A\|_*^\alpha}(R+\|x\|_*)^\alpha.
\end{align*}
This shows that $\mathcal{T}g$ is well defined and in particular $\mathcal T$ maps
$\mathcal C_{\alpha,R}$ into itself with
\begin{equation*}
 \|\mathcal Tg\|_{\alpha,R}\leq\kappa(C_{R}+L_f^*\|g\|_{\alpha,R}).
\end{equation*}
For $g,\widehat g\in\mathcal C_{\alpha,R}$, use
\[
 \|f(z+g(z))-f(z+\widehat g(z))\|_*
 \leq L_f^*\|g-\widehat g\|_{\alpha,R} (R+\|z\|_*)^\alpha
\]
in the two series.
Writing $w(x)=(R+\|x\|_*)^\alpha$ only in the following computation,
the stable component gives
\begin{align*}
 \| (\mathcal Tg-\mathcal T\widehat g)_s(x)\|_*
 &\leq L_f^*\|g-\widehat g\|_{\alpha,R}
       \sum_{j=1}^{\infty}\mathcal{C}_s^{j-1}w(A^{-j}x)\\
 &\leq L_f^*\|g-\widehat g\|_{\alpha,R}w(x)
       \sum_{j=1}^{\infty}\mathcal{C}_s^{j-1}\|A^{-1}\|_*^{\alpha j}\\
 &=L_f^*\|g-\widehat g\|_{\alpha,R}
       \frac{\|A^{-1}\|_*^\alpha}{1-\mathcal{C}_s\|A^{-1}\|_*^\alpha}w(x).
\end{align*}
For the unstable component, the corresponding computation is
\begin{align*}
 \| (\mathcal Tg-\mathcal T\widehat g)_u(x)\|_*
 &\leq L_f^*\|g-\widehat g\|_{\alpha,R}
       \sum_{j=0}^{\infty}\mathcal{C}_u^{j+1}w(A^jx)\\
 &\leq L_f^*\|g-\widehat g\|_{\alpha,R}w(x)
       \sum_{j=0}^{\infty}\mathcal{C}_u^{j+1}\|A\|_*^{\alpha j}\\
 &=L_f^*\|g-\widehat g\|_{\alpha,R}
       \frac{\mathcal{C}_u}{1-\mathcal{C}_u\|A\|_*^\alpha}w(x).
\end{align*}
Taking the maximum of the present components, dividing by $w(x)$,
and taking the supremum over $x$ gives the factor $\kappa$.
We obtain
\begin{equation}\label{eq:reverse-contraction}
 \|\mathcal Tg-\mathcal T\widehat g\|_{\alpha,R}
 \leq L_f^*\kappa\|g-\widehat g\|_{\alpha,R},
\end{equation}
which shows that $\mathcal{T}$ is contractive by the second inequality of \eqref{coeff.-contraction}. Consequently, $\mathcal{T}$ has a unique fixed point $g_*$ by the Banach fixed point theorem. 
Thus $G=\Id+g_*$ satisfies $G\circ A=F\circ G$.
\par
\medskip
\noindent\textit{Step 4. Verification that $H$ and $G$ are inverse maps.}
Lemmas~\ref{lem:composition-growth} and \ref{lem:equivalent-sublinear-growth}, together with \eqref{eq:explicit-growth-comparison},
\eqref{eq:h-estimate} and $g_*\in \mathcal{C}_{\alpha,R}$, allow us to write
\[
 H\circ G=\Id+p,\qquad G\circ H=\Id+k,\qquad
 \|p\|_{\alpha,R}<\infty,\quad\|k\|_{\alpha,R}<\infty.
\]
The two intertwining relations imply, with the compositions written out,
\[
 (H\circ G)\circ A
 =H\circ(G\circ A)
 =H\circ(F\circ G)
 =(H\circ F)\circ G
 =A\circ(H\circ G).
\]
Subtracting $A$ from both sides gives the relation for $p$:
\[
p\circ A=A \circ p.
\]
Consequently, for every $N\geq1$ and $x\in X$,
\[
 p_s(x)=A_s^Np_s(A^{-N}x),\qquad
 p_u(x)=A_u^{-N}p_u(A^Nx).
\]
By \eqref{eq:linear-weights}, for every $x\in X$,
\begin{align*}
 \|p_s(x)\|_*&\leq\|p\|_{\alpha,R}(\mathcal{C}_s\|A^{-1}\|_*^\alpha)^N(R+\|x\|_*)^\alpha,\\
 \|p_u(x)\|_*&\leq\|p\|_{\alpha,R}(\mathcal{C}_u\|A\|_*^\alpha)^N(R+\|x\|_*)^\alpha.
\end{align*}
Letting $N\to\infty$ and using \eqref{eq:norm-gaps}, we get $p=0$,
so $H\circ G=\Id$.
\par
For the other composition, the same relations give
\[
 (G\circ H)\circ F
 =G\circ(H\circ F)
 =G\circ(A\circ H)
 =(G\circ A)\circ H
 =F\circ(G\circ H).
\]
Since $G\circ H=\Id+k$ and $F=A+f$, this identity reads
\[
 Fx+k(Fx)=A(x+k(x))+f(x+k(x)) \qquad \text{for all } x\in X.
\]
Subtracting $Fx=Ax+f(x)$ yields the following equation.
Equivalently
\begin{equation}\label{eq:composition-error}
 k(Fx)=Ak(x)+r(x) \qquad \text{with } \qquad
 r(x)=f(x+k(x))-f(x).
\end{equation}
In particular,
\begin{equation}\label{eq:composition-error-bound}
 \|r(x)\|_*\leq L_f^*\|k(x)\|_*
 \leq L_f^*\|k\|_{\alpha,R} (R+\|x\|_*)^\alpha.
\end{equation}
Iterate \eqref{eq:composition-error} backward on $X_s$ and forward
on $X_u$.
We have 
\[
 k_s(x)=A_sk_s(F^{-1}x)+r_s(F^{-1}x),\qquad
 k_u(x)=A_u^{-1}k_u(Fx)-A_u^{-1}r_u(x).
\]
Inductively, for each $N\ge 1$ we derive 
\begin{align*}
 k_s(x)&=A_s^Nk_s(F^{-N}x)
          +\sum_{j=1}^{N}A_s^{j-1}r_s(F^{-j}x),\\
 k_u(x)&=A_u^{-N}k_u(F^Nx)
          -\sum_{j=0}^{N-1}A_u^{-j-1}r_u(F^jx).
\end{align*}
The boundary terms have the bounds
\begin{align*}
 \|A_s^Nk_s(F^{-N}x)\|_*
 &\leq\|k\|_{\alpha,R}(\mathcal{C}_sB^\alpha)^N(R+\|x\|_*)^\alpha,\\
 \|A_u^{-N}k_u(F^Nx)\|_*
 &\leq\|k\|_{\alpha,R}(\mathcal{C}_uD^\alpha)^N(R+\|x\|_*)^\alpha,
\end{align*}
and hence tend to zero by \eqref{eq:buffered-gaps}.
The remaining series converge by \eqref{eq:composition-error-bound}.
Let $N\to \infty$ implies that
\begin{align*}
 \|k_s(x)\|_*&\leq L_f^*\|k\|_{\alpha,R}
                 \frac{B^\alpha}{1-\mathcal{C}_sB^\alpha}(R+\|x\|_*)^\alpha,\\
 \|k_u(x)\|_*&\leq L_f^*\|k\|_{\alpha,R}
                 \frac{\mathcal{C}_u}{1-\mathcal{C}_uD^\alpha}(R+\|x\|_*)^\alpha.
\end{align*}
Taking the maximum, dividing by $(R+\|x\|_*)^\alpha$ and taking the supremum over $x$
yields
\[
 \|k\|_{\alpha,R}\leq L_f^*\kappa_{B,D}\|k\|_{\alpha,R}.
\]
Since $L_f^*\kappa_{B,D}<1$ given in \eqref{eq:buffered-smallness}, we conclude that $k=0$.
Therefore $G\circ H=H\circ G=\Id$. Both maps are continuous, so $H$
is a global homeomorphism with inverse $G$.
\par
\medskip
\noindent\textit{Step 5. Growth and uniqueness in the prescribed class.}
The bounds already obtained give
\begin{equation*}
 \|H(x)-x\|_*+\|H^{-1}(x)-x\|_*
 \leq C_{R}\left(\kappa_{B,D}
             +\frac{\kappa}{1-L_f^*\kappa}\right)(R+\|x\|_*)^\alpha.
\end{equation*}
Lemma \ref{lem:equivalent-sublinear-growth} implies \eqref{eq:two-sided-growth}, as desired. 

\par
Suppose that $\widehat H=\Id+\widehat h$ is another continuous solution
of \eqref{eq:semiconjugacy} such that, for some constant $C_{\widehat h}$,
\[
 \|\widehat h(x)\|\leq C_{\widehat h}(1+\|x\|^\alpha),\qquad x\in X.
\]
Then $d=\widehat H-H$ satisfies $d(Fx)=Ad(x)$ and $\|d\|_{\alpha,R}<\infty$.
For the latter assertion, norm equivalence first gives
\[
 \|\widehat h(x)\|_*
 \leq MC_{\widehat h}\max\{1,m^{-\alpha}\}
       (1+\|x\|_*^\alpha)
 \leq 2MC_{\widehat h}\max\{1,m^{-\alpha}\}
       (R+\|x\|_*)^\alpha.
\]
Together with \eqref{eq:h-estimate} this bounds the weighted norm of
$d=\widehat h-h$.
\par
Iterating its homogeneous equation yields for each $N\ge 1$
\[
 d_s(x)=A_s^Nd_s(F^{-N}x),\qquad
 d_u(x)=A_u^{-N}d_u(F^Nx).
\]
Exactly as for the boundary terms above, we obtain 
\begin{align*}
 \|d_s(x)\|_*&\leq\|d\|_{\alpha,R}(\mathcal{C}_sB^\alpha)^N(R+\|x\|_*)^\alpha,\\
 \|d_u(x)\|_*&\leq\|d\|_{\alpha,R}(\mathcal{C}_uD^\alpha)^N(R+\|x\|_*)^\alpha.
\end{align*}
Letting $N\to\infty$ proves $\widehat H=H$.
\par
Finally, the construction includes the pure contraction and expansion
cases by omitting the absent component. In the contraction case it also
gives an iteration-limit formula for the inverse conjugacy.
Indeed, $G\circ A=F\circ G$ gives
\[
 G(x)=F^n\bigl(A^{-n}x+g(A^{-n}x)\bigr).
\]
In this stable-only case, $L_f^*\kappa<1$ reads
$L_f^*<\|A^{-1}\|_*^{-\alpha}-\mathcal{C}_s$. Multiplication by
$\|A^{-1}\|_*^\alpha$ gives
$(\mathcal{C}_s+L_f^*)\|A^{-1}\|_*^\alpha<1$.
Since $\operatorname{Lip}_*(F)\leq \mathcal{C}_s+L_f^*$ and
$(\mathcal{C}_s+L_f^*)\|A^{-1}\|_*^\alpha<1$, one has
\begin{align*}
 \|G(x)-F^n(A^{-n}x)\|_*
 &\leq (\mathcal{C}_s+L_f^*)^n\|g(A^{-n}x)\|_*\\
 &\leq\|g\|_{\alpha,R}\bigl((\mathcal{C}_s+L_f^*)\|A^{-1}\|_*^\alpha\bigr)^n (R+\|x\|_*)^\alpha
 \longrightarrow0.
\end{align*}
Thus $G(x)=\lim_{n\to\infty}F^n(A^{-n}x)$, uniformly on bounded sets.
\par
In the pure expansion case, the smallness condition $L_f^*\kappa<1$ becomes
$L_f^*<\mathcal{C}_u^{-1}-\|A\|_*^\alpha$, so $\mathcal{C}_uL_f^*<1$ and
\[
 \|F(x)-F(y)\|_*\geq(\mathcal{C}_u^{-1}-L_f^*)\|x-y\|_*,\qquad
 \operatorname{Lip}_*(F^{-1})\leq\frac{\mathcal{C}_u}{1-\mathcal{C}_uL_f^*}.
\]
Here $L_f^*<\mathcal{C}_u^{-1}-\|A\|_*^\alpha$ also gives
\[
 \mathcal{C}_u\|A\|_*^\alpha<1-\mathcal{C}_uL_f^*,\qquad
 \frac{\mathcal{C}_u\|A\|_*^\alpha}{1-\mathcal{C}_uL_f^*}<1.
\]
Using $G(x)=F^{-n}(A^nx+g(A^nx))$ instead, we obtain
\[
 \|G(x)-F^{-n}(A^nx)\|_*
 \leq\|g\|_{\alpha,R}
       \left(\frac{\mathcal{C}_u\|A\|_*^\alpha}{1-\mathcal{C}_uL_f^*}\right)^n (R+\|x\|_*)^\alpha\longrightarrow0,
\]
and hence $G(x)=\lim_{n\to\infty}F^{-n}(A^nx)$, uniformly on bounded sets.
This completes the proof.
\end{proof}
\endgroup
\section*{Acknowledgements}
This work is supported by NSF-CQ \#CSTB2023NSCQ-JQX0020 and the Science and Technology Research Program of Chongqing Municipal Education Commission (Grant No. KJQN202300540).


\begin{thebibliography}{99}
\setlength{\baselineskip}{9.3pt}
\bibitem{AW06}B. Aulbach and T. Wanner, Topological simplification of nonautonomous difference equations,
{\it J. Difference Equ. Appl.} {\bf 12} (2006) 283--296.
\bibitem{BDP21}L. Backes, D. Dragi\v{c}evi\'c and K. J. Palmer, Linearization and H\"older continuity for nonautonomous systems,
{\it J. Differential Equations} {\bf 297} (2021) 536--574.
\bibitem{BV06} L. Barreira and C. Valls, A Grobman-Hartman theorem
for nonuniformly hyperbolic dynamics, {\it J. Differential Equations}, {\bf 228} (2006) 285-310.
\bibitem{BV06P}
L. Barreira and C. Valls, H\"older conjugacies for random dynamical systems,
{\it Phys. D} {\bf 223} (2006) 256-269.
\bibitem{BL94}P. W. Bates and K. Lu, A Hartman--Grobman theorem for the Cahn--Hilliard and phase-field equations,
{\it J. Dynam. Differential Equations} {\bf 6} (1994) 101--145.
\bibitem{B78}G. R. Belitskii, Equivalence and normal forms of germs of smooth mappings,
{\it Russian Math. Surveys} {\bf 33} (1978), no.~1, 107--177.
\bibitem{Brjuno7172}A. D. Brjuno, Analytic form of differential equations. I, II,
{\it Trudy Moskov. Mat. Ob\v{s}\v{c}.} {\bf 25} (1971) 119--262;
{\bf 26} (1972) 199--239 (in Russian).
\bibitem{BC06}X. Buff and A. Ch\'eritat, The Brjuno function continuously estimates the size of quadratic Siegel disks,
{\it Ann. of Math.} (2) {\bf 164} (2006) 265--312.
\bibitem{CR18}\'A. Casta\~neda and G. Robledo, Dichotomy spectrum and almost topological conjugacy on nonautonomus unbounded difference systems,
{\it Discrete Contin. Dyn. Syst.} {\bf 38} (2018) 2287--2304.
\bibitem{CMR07}E. A. Coayla-Teran, S.-E. A. Mohammed and P. R. C. Ruffino, Hartman--Grobman theorems along hyperbolic stationary trajectories,
{\it Discrete Contin. Dyn. Syst.} {\bf 17} (2007) 281--292.
\bibitem{DLLX25}D. Dragi\v{c}evi\'c, X. Li, W. Lu and Y. Xia, On the global linearization theorem of Charles Pugh,
{\it Proc. Roy. Soc. Edinburgh Sect. A}, First View (2025) 1--16.
\bibitem{E01}M. S. ElBialy, Local contractions of Banach spaces and spectral gap conditions,
{\it J. Funct. Anal.} {\bf 182} (2001) 108--150.
\bibitem{G59} D. M. Grobman, Homeomorphism of systems of differential equations (Russian),  {\it Dokl.
Akad. Nauk SSSR} {\bf 128} (1959) 880-881.
\bibitem{G62}D. M. Grobman, Topological classification of neighborhoods of a singularity in $n$-space,
{\it Mat. Sb. (N.S.)} {\bf 56(98)} (1962), no.~1, 77--94 (Russian).
\bibitem{H60} P. Hartman, A lemma in the theory of structural stability of differential equations
{\it Proc. Amer. Math. Soc.} {\bf 11} (1960) 610-620.
\bibitem{H60b}P. Hartman, On local homeomorphisms of Euclidean spaces,
{\it Bol. Soc. Mat. Mexicana (2)} {\bf 5} (1960) 220--241.
\bibitem{HP16}M.-L. Hein and J. Pr\"uss, The Hartman--Grobman theorem for semilinear hyperbolic evolution equations,
{\it J. Differential Equations} {\bf 261} (2016) 4709--4727.
\bibitem{Herman87}M. R. Herman, Recent results and some open questions on Siegel's linearization theorem of germs of complex analytic diffeomorphisms of $\mathbb C^n$ near a fixed point,
in {\it VIIIth International Congress on Mathematical Physics} (Marseille, 1986),
World Scientific, Singapore, 1987, pp. 138--184.
\bibitem{Huerta20}I. Huerta, Linearization of a nonautonomous unbounded system with nonuniform contraction: A spectral approach,
{\it Discrete Contin. Dyn. Syst.} {\bf 40} (2020) 5571--5590.
\bibitem{KP91}J. Kurzweil and G. Papaschinopoulos, Topological equivalence and structural stability for linear difference equations,
{\it J. Differential Equations} {\bf 89} (1991) 89--94.
\bibitem{Lin07}F. Lin, Hartman's linearization on nonautonomous unbounded system,
{\it Nonlinear Anal.} {\bf 66} (2007) 38--50.
\bibitem{Lu91}K. Lu, A Hartman--Grobman theorem for scalar reaction-diffusion equations,
{\it J. Differential Equations} {\bf 93} (1991) 364--394.
\bibitem{Palis68}J. Palis, On the local structure of hyperbolic points in Banach spaces,
{\it An. Acad. Brasil. Ci\^enc.} {\bf 40} (1968) 263--266.
\bibitem{Palmer73}K. J. Palmer, A generalization of Hartman's linearization theorem,
{\it J. Math. Anal. Appl.} {\bf 41} (1973) 753--758.
\bibitem{PM92}R. P\'erez-Marco, Solution compl\`ete au probl\`eme de Siegel de lin\'earisation d'une application holomorphe au voisinage d'un point fixe (d'apr\`es J.-C. Yoccoz),
{\it Ast\'erisque} {\bf 206} (1992), Exp.~753, 273--310.
\bibitem{Poincare1879}
H. Poincar\'e, {\it Sur les propri\'et\'es des fonctions d\'efinies
par les \'equations aux diff\'erences partielles},
th\`ese de doctorat, Facult\'e des Sciences de Paris,
Gauthier-Villars, Paris, 1879.
\bibitem{PR16}C. P\"otzsche and E. Russ, Topological decoupling and linearization of nonautonomous evolution equations,
{\it Discrete Contin. Dyn. Syst. Ser. S} {\bf 9} (2016) 1235--1268.
\bibitem{Pugh69}C. C. Pugh, On a theorem of P. Hartman,
{\it Amer. J. Math.} {\bf 91} (1969) 363--367.
\bibitem{RS04}H. M. Rodrigues and J. Sol\`a-Morales, Smooth linearization for a saddle on Banach spaces,
{\it J. Dynam. Differential Equations} {\bf 16} (2004) 767--793.
\bibitem{Siegel42}C. L. Siegel, Iteration of analytic functions,
{\it Ann. of Math.} (2) {\bf 43} (1942) 607--612.
\bibitem{Sternberg57} S. Sternberg, Local contractions and a theorem of Poincar\'e,
{\it Amer. J. Math.} {\bf 79} (1957) 809-824.
\bibitem{WX24}M. Wu and Y. Xia, Linearization of a nonautonomous unbounded system with hyperbolic linear part: A spectral approach,
{\it Discrete Contin. Dyn. Syst.} {\bf 44} (2024) 491--522.
\bibitem{WXZ25}M. Wu, Y. Xia and L. Zhou, Global linearization of nonautonomous difference equations in Banach space: Pugh-type theorem in unbounded settings,
SSRN preprint 5356076 (2025).
\bibitem{Yoccoz95}J.-C. Yoccoz, Th\'eor\`eme de Siegel, nombres de Bruno et polyn\^omes quadratiques,
{\it Ast\'erisque} {\bf 231} (1995) 1--88.
\bibitem{ZS20}
J. Zhao and J. Shen, A new proof for the Hartman-Grobman theorem for random dynamical systems,
{\it Proc. Amer. Math. Soc.} {\bf 148} (2020) 365-377.
\end{thebibliography}
\end{document}